\documentclass{svproc}
\makeatletter
\renewcommand*{\@fnsymbol}[1]{\ensuremath{\ifcase#1\or *\or \dagger\or \ddagger\or
   \mathsection\or \mathparagraph\or \|\or **\or \dagger\dagger
   \or \ddagger\ddagger \else\@ctrerr\fi}}
\makeatother

\usepackage{url}

\usepackage{amsfonts}

\usepackage{latexsym}
\usepackage{amssymb}
\usepackage{amsmath}
\usepackage[mathscr]{eucal}
\usepackage{graphicx}
\usepackage{hyperref}
\usepackage{caption}
\usepackage{subcaption}

\renewcommand{\qed}{\hfill{\ \ \rule{2mm}{2mm}} \vspace{0.2in}}
\newcommand{\ind}{1\hspace{-2.3mm}{1}}

\begin{document}
\mainmatter              
\title{Redundancy of Codes with Graph Constraints \thanks{An abridged version of this paper was presented in~\(23^{rd}-\)Thailand-Japan Conference on Discrete and Computational Geometry, Graphs, and Games (TJCDCGGG, 2021) and is available in the book of abstracts at \emph{https://www.math.science.cmu.ac.th/tjcdcggg/}}}
\titlerunning{Redundancy of Codes}

\author{Ghurumuruhan Ganesan \thanks{Corresponding Author}}
\authorrunning{G. Ganesan}

\institute{IISER, Bhopal-462066, \email{gganesan82@gmail.com}}

\maketitle
\date{}


\begin{abstract}
In this paper, we study the redundancy of linear codes with graph constraints.
First we consider linear parity check codes based on bipartite graphs with diversity
and with generalized graph constraints. We describe sufficient conditions on the constraint probabilities and use the probabilistic
method to obtain linear codes that achieve the Gilbert-Varshamov redundancy bound in addition to
satisfying the constraints and the diversity index. In the second part we consider a generalization of graph capacity which we call as the fractional graph capacity and use the probabilistic method to determine bounds on the fractional capacity for arbitrary graphs. Specifically, we establish an upper bound in terms of the full graph capacity and a lower bound in terms of the average and maximum vertex degree of the graph.


\vspace{0.1in} \noindent \textbf{Key words:}  Linear codes, bipartite graphs, fractional graph capacity.

\end{abstract}

\bigskip

\setcounter{equation}{0}
\renewcommand\theequation{\arabic{section}.\arabic{equation}}
\section{Introduction} \label{intro}
Codes based on graphs arise often in both theory and applications and it is important to understand redundancies of such codes. Typically, redundancy bounds like Gilbert-Varshamov, Hamming and Singleton are obtained under the Hamming distance measure with no restrictions on the codes themselves. In many applications, the code itself might have additional graph constraints.

In this paper we study redundancy of codes with graph constraints as described in the following two subsections. 




\subsection*{Linear Parity Check Codes}
In the first part of the paper, we study redundancy of linear parity check codes based on bipartite graphs with diversiy and with general constraints. Linear codes with graph based constraints arise often in applications. For example, low density parity check (LDPC) codes~\cite{urb} that are used extensively in communication systems have the constraint that the left and right vertex degrees in the bipartite graph representation are small compared to the total number of vertices. Another example is expander codes~\cite{sipser} that requires that the bipartite graph representation to satisfy an \emph{expansion} property with respect to the left nodes. A linear time encoding and decoding algorithm for expander codes is also presented in~\cite{sipser}. Recently,~\cite{hemen} has proposed a local-decoding algorithm for expander codes capable of correcting a constant fraction of errors using a (relatively) small number of symbols from corrupted codeword.

In Section~\ref{sec_lin_code_gr}, we consider linear codes with diversity and generalized graphical constraints and use random bipartite graphs to show that if the constraints are not severe with regards to their probability of occurrence, then there are codes with diversity that attain the Gilbert-Varshamov redundancy bound in addition to satisfying the said constraints.

\subsection*{Fractional Graph Capacity}
The capacity of a graph was introduced in~\cite{shannon} to determine the limits of error-free communication across channels with an underlying confusion graph. Lov\'asz~\cite{lazlo} obtained the expression for capacity of the cycle on~\(5\) vertices using ``umbrella" projection methods and since then, many bounds as well as variants of the capacity have been studied. For example,~\cite{hell} studied analogues of Shannon capacity and its connections with the ultimate chromatic number. Marton~\cite{marton} obtained expressions for graph capacities for a sequence of graphs based on typical sequences and more recently,~\cite{alon} investigated the problem of approximating graph capacity by finite graph products.


In this paper, we introduce the notion of fractional capacity of a graph and obtain bounds in terms of the graph degree parameters. Physically, fractional capacity corresponds to the communication limits of channels that corrupt at most a \emph{fraction} of symbols in a long codeword. We first use the probabilistic method to obtain upper bound for the fractional capacity in terms of the full graph capacity and a lower bound in terms of the average and maximum vertex degrees. We then consider a special class of graphs called square graphs and determine an upper bound for the fractional capacity in terms of the minimum degree of its root. In particular for regular square graphs, this obtains upper and lower bounds for the fractional capacity in terms of the common vertex degree.

The paper is organized as follows: In Section~\ref{sec_lin_code_gr} we study linear parity check codes with graphical constraints and use random bipartite graphs to determine achievability of the Gilbert-Varshamov bound. Next in Section~\ref{sec_frac_cap}, we define fractional capacity of a graph and state and prove our main result involving upper and lower bounds for the fractional capacity.


\setcounter{equation}{0}
\renewcommand\theequation{\arabic{section}.\arabic{equation}}
\section{Linear Parity Check Codes with Graph Constraints}\label{sec_lin_code_gr}
We begin with some general definitions. For a finite set~\({\cal X}\) and integer~\(n \geq 1,\) an~\(n-\)length word (based on~\({\cal X}\)) is an element of~\({\cal X}^{n}.\) An~\(n-\)length code~\({\cal C}\) is a subset of~\({\cal X}^{n}.\) If~\({\cal X}\) is a finite field, then we have the concept of linear codes: we say that~\({\cal C}\) is \emph{linear} if for any two words~\(\mathbf{c},\mathbf{d} \in {\cal C}\) and any two scalars~\(\alpha,\beta \in {\cal X},\) the word~\(\alpha \cdot \mathbf{c} + \beta \cdot \mathbf{d} \in {\cal C}\) as well.

For two words~\(\mathbf{c} = (c_1,\ldots,c_n)\) and~\(\mathbf{d} = (d_1,\ldots,d_n)\) in~\({\cal X}^{n},\) we define the \emph{Hamming distance} between~\(\mathbf{c}\) and~\(\mathbf{d}\) to be
\begin{equation}\label{dist_def}
d(\mathbf{c},\mathbf{d}) = \sum_{i=1}^{n} \ind(c_i \neq d_i),
\end{equation}
where~\(\ind(.)\) refers to the indicator function. In this section, all distances are Hamming. We define the minimum distance~\(d_{H}({\cal C})\) of the code~\({\cal C}\) to be the minimum distance between any two codewords of~\({\cal C}\) and set the relative distance of~\({\cal C}\) to be
\begin{equation}\label{rel_dist_uni}
\delta_H({\cal C}) := \frac{d_{H}({\cal C})-1}{n}.
\end{equation}
The rate and redundancy of~\({\cal C}\) are respectively defined as
\begin{equation}\label{rate_def_uni}
R({\cal C}) := \frac{\log(\#\cal C)}{n\log(\#{\cal X})}
\end{equation}
and
\begin{equation}\label{red_def_uni}
\xi({\cal C}) := 1-R({\cal C}) = 1-\frac{\log(\#\cal C)}{n\log(\#{\cal X})},
\end{equation}
where logarithms are to the base two throughout and~\(\#{\cal C}\) denotes the size of~\({\cal C}.\)


In this section we consider binary linear codes with~\({\cal X} = \{0,1\}\) and begin with a description of the random graph construction of linear parity check codes. Consider a random bipartite graph with left vertex set~\(X = \{a_1,a_2,\ldots,a_n\}\) and right vertex set~\(Y = \{b_1,b_2,\ldots,b_m\}\) obtained as follows.  Let~\(\{X_{i,j}\}_{1 \leq i \leq n, 1 \leq j \leq m}\) be independent and identically distributed binary random variables with
\[\mathbb{P}(X_{1,1} = 1) = p =  1-\mathbb{P}(X_{1,1} = 0)\]
where~\(0 < p < \frac{1}{2}\) is a constant. Throughout, constants do not depend on~\(n.\)
An edge is present between vertices~\(a_i\) and~\(b_j\) if and only if~\(X_{i,j}=1.\)
Let~\(m=n\epsilon\) for some constant~\(\epsilon  >0\) and let~\(G = G_n\) be the resulting random graph defined on the probability space~\((\Omega_n,{\cal F}_n, \mathbb{P}_n).\) For simplicity we drop the subscript from~\(\mathbb{P}_n\) henceforth.

For~\(1 \leq j \leq m\) let~\({\cal R}_j\) be the set of neighbours of the right vertex~\(b_j\)
and define the (random) code~\({\cal C}\) as follows. A word~\(\mathbf{c} = (c_1,\ldots,c_n) \in {\cal C}\) if and only if
\begin{equation}\label{code_def}
\oplus_{i \in {\cal R}_j} c_i = 0 \text { for all }1 \leq j \leq m.
\end{equation}
By construction (see~\cite{sipser}), the code~\({\cal C}\) is linear with rate at least~\(1-\frac{m}{n} =1-\epsilon.\)

We now introduce two restrictions on~\({\cal C}.\) The first restriction which we call as diversity is an expansion-type property with regards to the right nodes of the bipartite graph. The other restriction which we call as constraints are general events pertaining to the random graph as constructed in the previous paragraph.

For~\(0 < \gamma < 1\) we say that~\({\cal C}\) has a \emph{diversity index} of at least~\(\gamma\) if
\begin{equation}\label{exp_rev}
\#\left({\cal R}_x \setminus {\cal R}_y \right) \geq \gamma \#{\cal R}_x \text{ for any } b_x,b_y \in Y.
\end{equation}
Thus any two parity nodes have at least a fraction~\(\gamma\) of different neighbours. We could think of condition~(\ref{exp_rev}) as a mild form of the expansion property with respect to the \emph{parity nodes}. We remark that in the usual construction via expander graphs, the condition for expansion is with regards to the (left) codeword index nodes of the bipartite graphs (see~\cite{sipser}).

We now define constraints on linear codes as follows. A~\(n-\)length \emph{constraint}~\({\cal E}_n\) is an event in~\({\cal F}_n.\) For example, the event~\({\cal H}_n\) that for each~\(2 \leq i \leq \sqrt{n}\) the right vertices~\(b_{i-1}\) and~\(b_{i+1}\) both have the left vertex~\(a_i\) as a neighbour is an example of a constraint. We say that the random graph~\(G\) satisfies the constraint~\({\cal E}_n\) if~\(G \in {\cal E}_n.\) Given a collection of constraints~\({\cal E}_n\) and real numbers~\(0 < \delta, \gamma < 1,\) we would like to know if there is a linear code with relative distance~\(\delta\) and diversity index~\(\gamma\) that also satisfies the constraints. If so, what would be the redundancy of such a code?

If there were no constraints or diversity, then the Gilbert-Varshamov bound (Theorem~\(4.2.1,\)~\cite{guru}) implies that there exists an~\(n-\)length code with redundancy at most~\(H(\delta)+o(1),\) where~\(o(1) \longrightarrow 0\) as~\(n \rightarrow \infty\) and
\begin{equation}\label{ent_def}
H(x) := -x \cdot \log{x} - (1-x) \cdot \log(1-x)
\end{equation}
is the (binary) entropy function. Throughout logarithms are to the base~\(2.\) Does imposing diversity and constraints increase the redundancy of a linear code? The following result says that if the constraints are not too strict, then we can still get linear parity check codes satisfying the Gilbert-Varshamov redundancy bound and with a given diversity index. Constants mentioned throughout do not depend on~\(n.\)
\begin{theorem}\label{main_thm} Let~\(0 <\delta< \frac{1}{2}\)  and~\(0 < \gamma < 1\) be any two constants and let~\(\{{\cal E}_n\}\) be a collection of constraints such that the probability~\(p_n = \mathbb{P}({\cal E}_n)\) satisfies
\begin{equation}\label{cons_bd}
\frac{\log\left(\frac{1}{p_n}\right)}{n} \longrightarrow 0
\end{equation}
as~\(n\rightarrow \infty.\)

There exists a deterministic~\(n-\)length linear parity check code~\({\cal D}_n\) with relative distance at least~\(\delta,\) diversity index at least~\(\gamma,\) redundancy~\(H(\delta)+o(1)\) and satisfying the constraint~\({\cal E}_n.\)
\end{theorem}
The condition~(\ref{cons_bd}) is satisfied, for example, if~\[\mathbb{P}(E_n) \geq e^{-f(n)}\] for some function~\(f\) such that~\(\frac{f(n)}{n} \longrightarrow 0\) as~\(n \rightarrow \infty.\) For all~\(n\) large, the code~\({\cal D}_n\) then attains the Gilbert-Varshamov bound in addition to satisfying the constraint~\({\cal E}_n.\)

For example, the event~\({\cal H}_n\) described prior to the statement of Theorem~\ref{main_thm} occurs with probability at least~\(p^{2\sqrt{n}}\) and so satisfies~(\ref{cons_bd}). Consequently, for all~\(n\) large, there exists a linear code with relative distance at least~\(\delta,\) diversity index at least~\(\gamma,\) redundancy~\(H(\delta)+o(1)\) and satisfying the constraint~\({\cal H}_n.\) 

Below, we obtain the desired code in Theorem~\ref{main_thm} by the probabilistic method.
\subsection*{Proof of Theorem~\ref{main_thm}}\label{pf_main_thm}

The proof of Theorem~\ref{main_thm} consists of three steps. In the first step, we choose the edge probability~\(p\) to be an appropriate constant so that the diversity condition is ensured. For such a choice of~\(p,\) we show in the second step that the minimum distance of the code~\({\cal C}\) as obtained in~(\ref{code_def}) is at least~\(\delta n+1\) with high probability i.e., with probability converging to one as~\(n \rightarrow \infty.\) Finally, in the third step, we incorporate the constraints into~\({\cal C}.\) Throughout we let~\(\epsilon > H(\delta)\) be a constant and let~\(m=n\epsilon\) be the number of parity (right) nodes in the graph~\(G\) so that the size of~\({\cal C}\) is at least~\(2^{n(1-\epsilon)}.\)

\emph{\underline{Step 1 (Ensuring diversity)}}: Let~\({\cal C}\) be the linear code as obtained in~(\ref{code_def}). To ensure that~\({\cal C}\) satisfies the diversity property, we argue as follows. Let~\(b_x\) and~\(b_y\) be any two right vertex nodes. We have that a left vertex~\(a_i\) is present in~\({\cal R}_x\) with probability~\(p\) and is present in~\({\cal R}_x \cap {\cal R}_y\)
with probability~\(p^2.\) Therefore by standard deviation estimates (Corollary~\(A.1.14,\) pp.~\(312,\)~\cite{alon2}),
we have for~\(0 < \theta < \frac{1}{4}\) that
\begin{equation}\label{eqr_one}
\mathbb{P}\left(\left|\#{\cal R}_x -np\right| \geq np\theta\right) \leq \exp\left(-\frac{\theta^2}{4}np\right)
\end{equation}
and that
\begin{equation}\label{eqr_two}
\mathbb{P}\left(\left|\#({\cal R}_x \cap {\cal R}_y) -np^2\right| \geq np^2\theta\right) \leq \exp\left(-\frac{\theta^2}{4}np^2\right).
\end{equation}
Letting
\[R_{tot} := \bigcap_{x} \{\left|\#{\cal R}_x -np\right| \geq np\theta\}\] we then get that
\begin{equation}\label{r_tot_est}
\mathbb{P}(R_{tot}) \geq 1-2m^2e^{-\frac{\theta^2}{4}np^2}.
\end{equation}

Similarly, using~\(np^2 < np,\) we get from~(\ref{eqr_one}) and~(\ref{eqr_two}) that the event
\[F_{x,y} := \left\{\#\left({\cal R}_x \setminus {\cal R}_y \right) \geq np(1-\theta) - np^2(1+\theta)\right\}\]
occurs with probability at least~\(1-2e^{-\frac{\theta^2}{4}np^2}\)
and so letting~\(F_{tot} := \bigcap_{x,y} F_{x,y},\) we then get that
\begin{equation}\label{f_tot_est}
\mathbb{P}(F_{tot}) \geq 1-2m^2e^{-\frac{\theta^2}{4}np^2}.
\end{equation}
From~(\ref{r_tot_est}),~(\ref{f_tot_est}) and the union bound, we therefore we get that the event\\\(E_{div} := R_{tot} \cap F_{tot}\) occurs with probability
\begin{equation}\label{e_div_est}
\mathbb{P}(E_{div}) \geq 1-4m^2e^{-\frac{\theta^2}{4}np^2}.
\end{equation}

If~\(E_{div}\) occurs, then for any right nodes~\(b_x,b_y\)  we have
\begin{eqnarray}
\frac{\#\left({\cal R}_x \setminus {\cal R}_y \right)}{\#{\cal R}_x} &\geq& \frac{np(1-\theta) - np^2(1+\theta)}{np(1+\theta)} \nonumber\\
&=& \frac{1-\theta}{1+\theta} - p \nonumber
\end{eqnarray}
which is at least~\(\gamma\) provided~\(\theta,p\) are sufficiently small constants. We henceforth fix such a~\(p.\)

\emph{\underline{Step 2 (Estimating the minimum distance)}}: For a set~\({\cal S} \subseteq \{1,2,\ldots,n\}\)
we define the word~\(\mathbf{v}({\cal S}) = (v_1,\ldots,v_n)\) satisfying
\begin{equation}\label{vs}
v_i  = \left\{
\begin{array}{cc}
1 & \text{ if } i \in {\cal S}\\
0 & \text{ otherwise}.
\end{array}
\right.
\end{equation}

Letting~\({\cal S} = \{a_{l_1},\ldots,a_{l_g}\}\) be any set of left vertices
we upper bound the probability that~\({\cal S}\) would cause no parity check violations; i.e.
we estimate the probability that the word~\(\mathbf{v}({\cal S})\) (see~(\ref{vs})) belong to the code~\({\cal C}.\)
We consider two cases depending on whether~\(\#{\cal S} \leq t\) or not, for some integer
constant~\(t \geq 1\) to be determined later.

\emph{Case I~(\(\#{\cal S}  = g \leq t\)):} For~\(1 \leq i \leq g\) let~\({\cal N}_i\) be the set of neighbours of the left vertex~\(a_{l_i}.\) Let~\({\cal M}_i := {\cal N}_{i} \setminus \left(\bigcup_{\stackrel{1 \leq \l \leq t}{l\neq i}} {\cal N}_{l}\right)\) be the set
of (unique) neighbours of vertex~\(a_{l_i}\) not adjacent to any of the remaining vertices in~\({\cal S} \setminus \{a_{l_i}\}.\)
Defining~\(Y_{i,j} := \ind\left(b_j \in {\cal M}_i\right)\) we then have that~\(\{Y_{i,j}\}_{1 \leq j \leq m}\)
are independent and identically distributed for any~\(1 \leq i \leq g,\) with
\[\mathbb{P}(Y_{i,j} = 1) = p(1-p)^{g-1} = 1-\mathbb{P}(Y_{i,j} = 0).\]
Thus~\(\mathbb{P}\left({\cal M}_i = \emptyset\right) = \left(1-p(1-p)^{g-1}\right)^{m} \leq e^{-mp(1-p)^{g-1}}\)
and consequently
\begin{equation}\label{mi_est}
\mathbb{P}\left(\bigcup_{1 \leq i \leq g} \{{\cal M}_i = \emptyset\}\right) \leq g e^{-mp(1-p)^{g-1}}.
\end{equation}


If the event~\(E\left({\cal S}\right) := \bigcap_{1 \leq i \leq g} \{{\cal M}_i \neq \emptyset\}\) occurs,
then the word~\(\mathbf{v}({\cal S}) \notin {\cal C}.\) Therefore if the event
\begin{equation}\label{e_low_def}
E_{low} := \bigcap_{{\cal S}} E\left({\cal S}\right)
\end{equation}
occurs where the intersection is with respect to all subsets~\({\cal S} \subset \{1,2,\ldots,n\}\)
of size~\(g \leq t,\) then the minimum distance of any word in~\({\cal C}\) from the all zeros codeword is
at least~\(t+1.\) Since~\({\cal C}\) is linear this implies that the minimum distance of~\({\cal C}\)
is at least~\(t+1.\)

We now see that~\(E_{low}\) occurs with high probability, i.e., with probability
converging to one as~\(n \rightarrow \infty.\)
If~\(E_{low}^c\) denotes the complement of the set~\(E_{low},\)
then from~(\ref{mi_est}) we have
\begin{equation}\nonumber
\mathbb{P}\left(E_{low}^c\right) \leq \sum_{g=1}^{t} g {n \choose g} e^{-mp(1-p)^{g-1}}\leq t^2 {n \choose t} e^{-mp(1-p)^{t-1}}
\end{equation}
provided~\(t < \frac{n}{2}.\) Using~\({n \choose t} \leq \left(\frac{ne}{t}\right)^{t}\) we further get
that~
\begin{equation}\label{e_est2}
\mathbb{P}(E_{low}^c) \leq e^{-\Delta_0}
\end{equation}
where~\[\Delta_0 := mp(1-p)^{t-1} - t \log\left(\frac{ne}{t}\right) - 2\log{t}.\]
Since~\(m = \epsilon n\) and~\(p >0\) is a constant, we have that
\begin{equation}\label{del_est}
\Delta_0 \geq \frac{m}{2}p(1-p)^{t-1} \geq 4C \cdot n
\end{equation}
for all~\(n\) large and some constant~\(C > 0.\)

\emph{Case II~(\( t+1 \leq \#{\cal S}  \leq \delta n\)):}
For a right vertex~\(b_j,\) we recall that~\({\cal R}_j\)
is the random set of (left) neighbours of the vertex~\(b_j.\)
Define the event
\[F_j({\cal S}) := \{\#\left({\cal R}_j \cap {\cal S}\right) \text{ is odd}\}.\]

If~\(F_j({\cal S})\) occurs, then the word~\(\mathbf{v}({\cal S})\) would cause
a parity check violation at the right vertex~\(b_j.\)
Therefore if~\(\bigcup_{1 \leq j \leq m} F_j({\cal S})\) occurs,
then~\(\mathbf{v}({\cal S}) \notin {\cal C}.\)
Define the event
\begin{equation}\label{e_up_def}
E_{up} := \bigcap_{{\cal S}} \left(\bigcup_{1 \leq j \leq m} F_j({\cal S})\right)
\end{equation}
where the intersection is with respect to all sets~\({\cal S}\)
whose cardinality lies between~\(t+1\) and~\( \delta n.\)
Extending the above argument we see that if~\(E_{up}\) occurs, then there is
no word in~\({\cal C}\) whose distance from the all zeros codeword lies between~\(t+1\) and~\(\delta n.\)
Combining with the event~\(E_{low}\) defined in~(\ref{e_low_def}),
we have that if~\(E_{low} \cap E_{up}\) occurs,
then the minimum distance of the code is at least~\(\delta n + 1.\)

We estimate the probability that~\(E_{up}\) occurs.
For any right vertex~\(b_j,\) the number of left neighbours~\(\#{\cal R}_j\)
is Binomially distributed with parameters~\(n\) and~\(p.\)
Therefore for a deterministic set~\({\cal S}\) with~\(\#{\cal S} = g,\)
the cardinality of the random set~\(\#({\cal R}_j \cap {\cal S})\) is
Binomially distributed with parameters~\(g\) and~\(p.\)
Therefore
\begin{eqnarray}
\mathbb{P}(F_j^c({\cal S})) &=& \sum_{\stackrel{0 \leq k \leq g}{k \text{ even }}} {g \choose k} p^{k}(1-p)^{g-k} \nonumber\\
&=& \frac{1}{2}\left((1-p+p)^{g} + ((1-p)-p)^{g}\right) \nonumber\\
&=&\frac{1}{2} + \frac{1}{2}(1-2p)^{g}. \label{fj_comp}
\end{eqnarray}

Let~\(0 < \eta < \frac{1}{2}\) be a small constant.  Using~\(g \geq t+1\) and choosing~\(t\) sufficiently large, we then get
from~(\ref{fj_comp}) that~\[\mathbb{P}\left(F_{j}^c({\cal S})\right) \leq \frac{1}{2^{1-\eta}}\]
and so~\[\mathbb{P}\left(\bigcap_{j=1}^{m} F_{j}^c({\cal S})\right) \leq \left(\frac{1}{2^{1-\eta}}\right)^{m} = \frac{1}{2^{(1-\eta)\epsilon n}}.\]
There are~\({n \choose g}\) sets of cardinality~\(g\) and so from~(\ref{e_up_def}), we therefore have that
\begin{equation}\label{e_up_est2}
\mathbb{P}\left(E_{up}^c\right) \leq \left(\sum_{g=t+1}^{\delta n}{n \choose g}\right) \cdot \frac{1}{2^{(1-\eta)\epsilon n}}\leq
\frac{1}{2^{\beta n}},
\end{equation}
where~\(\beta := (1-\eta)\epsilon - H(\delta)\) and the final inequality in~(\ref{e_up_est2}) follows from  standard Hamming ball estimates (Proposition~\(3.3.1\)~\cite{guru}).  Since~\(\epsilon > H(\delta)\) strictly, we choose~\(\eta > 0\) small enough so that~\(\beta\)
is strictly positive. Fixing such an~\(\eta,\) we define~\(E_{dist} := E_{low} \cap E_{up}\) and have from~(\ref{e_up_est2}),~(\ref{e_est2}) and~(\ref{del_est}) that
\begin{equation}\label{e_low_up}
\mathbb{P}(E_{dist}) \geq 1-e^{-4Cn} - \frac{1}{2^{\beta n}} \geq 1-e^{-3Cn}
\end{equation}
for all~\(n\) large, where the constant~\(C >0\) is as in~(\ref{del_est}).

Combining~(\ref{e_div_est}) and~(\ref{e_low_up}) and using the fact that~\(m = \epsilon n,\) we get from a union bound that
\begin{equation}\label{e_div_low}
\mathbb{P}(E_{div} \cap E_{dist}) \geq 1-4m^2e^{-\frac{\theta^2}{4}np^2} - e^{-3Cn} \geq 1-e^{-2Dn}
\end{equation}
for all~\(n\) large and some constant~\(D > 0.\)

\emph{\underline{Step 3 (Incorporating constraints)}}: From~(\ref{cons_bd}), we have that~\({\cal E}_n\) occurs with probability at least~\(e^{-Dn}\)  for all~\(n\) large, where~\(D > 0\) is as in~(\ref{e_div_low}). Therefore from~(\ref{e_div_low}) we have~\[\mathbb{P}\left({\cal E}_n \cap E_{div} \cap E_{dist}\right) \geq e^{-Dn} - e^{-2Dn} > 0\]
and this implies that there exists an~\(n-\)length linear code with relative distance at least~\(\delta,\) diversity index at least~\(\gamma,\)
redundancy at most~\(\epsilon,\) and satisfying the constraint~\({\cal E}_n.\)~\(\qed\)

\setcounter{equation}{0}
\section{Fractional Graph Capacity}\label{sec_frac_cap}
Let~\(G = (V,E)\) be any connected graph containing~\(\#V = n \geq 3\) vertices and for a vertex~\(u,\) let~\({\mathcal N}_{G}[u]\)
be the set of all neighbours of~\(u,\) including~\(u.\) The graph distance between any two nodes~\(u\) and~\(v\) is the number of edges in the shortest
path between~\(u\) and~\(v.\) A set~\({\mathcal F} \subset V\) of vertices is said to be stable if no two vertices in~\({\mathcal F}\) are adjacent to each other in~\(G.\) We denote~\(\alpha(G)\) to be the independence number, i.e. the maximum size of a stable set in~\(G.\)

For integer~\(r \geq 1\) let~\(G(r)\) be the~\(r^{th}\) strong graph product of~\(G\) obtained as follows: The graph~\(G(r)\) has vertex set~\(V^{r}\)
and two vertices~\(\mathbf{u}= (u_1,\ldots,u_r)\) and~\(\mathbf{v} = (v_1,\ldots,v_r)\) are adjacent
if and only if~\(u_i \in {\mathcal N}_G[v_i]\) for each~\(1 \leq i \leq r.\)  For an integer~\(1 \leq k \leq r,\) we now define the subgraph~\(G(r,k) \subseteq G(r)\) as follows. Two vertices~\(\mathbf{v}, \mathbf{u} \in V^r\) are adjacent in~\(G(r,k)\) if and only if~\(\mathbf{u}\) and~\(\mathbf{v}\) are adjacent in~\(G(r)\) and differ in at most~\(k\) entries i.e.,~\[\sum_{i=1}^{r} \ind(u_i \neq v_i) \leq k,\] where~\(\ind(.)\) refers to the indicator function. We have the following definition.
\begin{definition}\label{frac_graph}
For a real number~\(0 < \gamma \leq 1\) we define the~\(\gamma-\)fractional capacity of~\(G\) to be
\begin{equation}\label{del_frac}
\Theta_{\gamma}(G) := \sup_{r \geq \frac{1}{\gamma}} \left(\alpha(G(r,\gamma r))\right)^{\frac{1}{r}}.
\end{equation}
\end{definition}
For~\(\gamma =1,\) the term~\(\Theta_1(G) =: \Theta(G)\) is the graph capacity  as defined in~\cite{shannon}. For differentiation, we refer to~\(\Theta(G)\) as the \emph{full} graph capacity and~\(\Theta_{\gamma}(G)\) as the \emph{fractional}  graph capacity.

In the context of codes, each vertex of~\(G(r)\) is a codeword of length~\(r\) and the term~\(\gamma\) represents the maximum fraction of symbols  that undergo corruption when passed through a channel with confusion graph~\(G.\) The quantity~\((\Theta_{\gamma}(G))^{r}\) is then the maximum size of a code from~\(G(r)\) that allows for error free communication.




For~\(0 \leq x \leq  1\) we let~\(H(x)\) be the entropy function as in~(\ref{ent_def}) and have the following result.
\begin{theorem}\label{main_thm0} Let~\(G\) be a connected graph on~\(n\) vertices and let~\(d_{av}\) and~\(\Delta\) be the average and maximum vertex degree of~\(G,\) respectively. For any~\(0 < \gamma \leq 1\) we have that
\begin{equation}\label{thet_bds}
n \cdot \max(f(\gamma,d_{av}),f(\gamma,\Delta),f(\gamma,n-1)) \leq \Theta_{\gamma}(G) \leq n \cdot \left(\frac{\Theta(G)}{n}\right)^{\gamma}
\end{equation}
where
\begin{equation}\label{low_bound}
f(\gamma,x) :=
\left\{
\begin{array}{cc}
\left(2^{H(\gamma)} \cdot x^{\gamma} \right)^{-1} & \text{ for } 0 < \gamma  < \frac{x}{x+1}\\
&\\
(x+1)^{-1} & \text{ for } \frac{x}{x+1} \leq \gamma \leq 1.
\end{array}
\right.
\end{equation}
\end{theorem}
We have the following remarks:\\
\emph{\underline{Remark 1}}: From~(\ref{thet_bds}) and the fact that~\(d_{av} \geq 1\) (recall that~\(G\) is connected) we see that for~\(0 < \gamma \leq \frac{1}{2},\) the fractional graph capacity grows at least as~\(\frac{n}{2^{H(\gamma)} \cdot d_{av}^{\gamma}}.\) For example the graph~\(G = C_5,\) the cycle on~\(5\) vertices, has~\(d_{av} = \Delta = 2\)  and it is well-known~\cite{lazlo} that the full graph capacity~\(\Theta(G) = \sqrt{5}\). Therefore setting~\(\gamma = \frac{1}{2},\) we get from~(\ref{thet_bds}) that the half graph capacity of~\(C_5\) satisfies
\[ \frac{5}{2\sqrt{2}} \leq \Theta_{\frac{1}{2}}(C_5) \leq \frac{5}{\sqrt[4]{5}}.\]

\emph{\underline{Remark 2}}: In general, we see from~(\ref{thet_bds}) that as~\(\gamma \rightarrow 0,\) the fractional graph capacity~\(\Theta_{\gamma}(G) \rightarrow n,\) the maximum possible value. On the other end, setting~\(\gamma =1 \) in the lower bound~(\ref{low_bound}), we get that the full graph capacity~\[\Theta(G) \geq \frac{n}{d_{av}+1}\] which is also obtained via the Tur\'an's bound~\cite{west}.

\emph{\underline{Remark 3}}: From~(\ref{del_frac}), we see that the upper bound for the fractional graph capacity is in terms of the full graph capacity, while the lower bound is in terms of vertex degrees. In the next section, we consider a particular class of graphs called square graphs and obtain upper \emph{and} lower bounds for the fractional capacity in terms of the vertex degrees.

\subsection*{Proof of Theorem~\ref{main_thm0}}
We prove the lower bound in~(\ref{thet_bds}) using the probabilistic method and a maximal stable set argument (the Gilbert-Varshamov argument~\cite{huff}) and prove the upper bound in~(\ref{thet_bds}) using a recursion estimate similar to the Singleton argument~\cite{huff}.\\
\emph{Proof of the lower bound in~(\ref{thet_bds})}: For a uniformly random vector~\(\mathbf{v}  = (v_1,\ldots,v_r) \in V^{r},\) let~\({\mathcal B}_{\gamma}(\mathbf{v})\) be the set of all vertices adjacent to~\(\mathbf{v}\) in the graph~\(G(r,\gamma r).\) The vertices~\(\{v_{j}\}_{1 \leq j \leq n}\) are mutually independent and uniformly distributed in~\(V\) and so the expected degree of~\(v_{j}\) is~\(d_{av}.\) Consequently, the expected size of~\({\mathcal B}_{\gamma}(\mathbf{v})\) is
\begin{equation}\label{eq_one2}
\mathbb{E}\#{\mathcal B}_{\gamma}(\mathbf{v}) = \sum_{k=0}^{\gamma r} {r \choose k} d^{k}_{av}.
\end{equation}

For~\(0 < \gamma < 1-\frac{1}{d_{av}+1},\) we use the Hamming ball estimate (see Proposition~3.3.1~\cite{guru}) to get that
\begin{equation}\label{eq_two2}
\mathbb{E}\#{\mathcal B}_{\gamma}(\mathbf{v}) = \sum_{k=0}^{\gamma r}{r \choose k} d_{av}^{k} \leq n^{\theta r}
\end{equation}
where~\(\theta = \frac{H(\gamma) + \gamma \log{d_{av}}}{\log{n}} \)
satisfies~\(0 < \theta  <1.\) For~\(0 < \epsilon = \epsilon(r) < 1\) to be determined later, we let
\begin{equation}\label{a_set_def}
{\mathcal A}(\epsilon) := \{ \mathbf{u} \in V^{r} : \#{\mathcal B}_{\gamma}(\mathbf{u}) \leq n^{r(\theta+\epsilon)}\}
\end{equation}
and obtain from Markov inequality and~(\ref{eq_two2}) that
\begin{equation}\label{a_eps_est}
\#{\mathcal A}(\epsilon) \geq n^{r}\left(1 - n^{-r\epsilon}\right).
\end{equation}

Let~\({\mathcal D} := \{\mathbf{w}_1,\ldots,\mathbf{w}_M\} \subseteq {\mathcal A}(\epsilon)\) be a stable set of maximum size in~\(G(r,\gamma r).\) By the maximality, we must have that the union~\(\bigcup_{i=1}^{M} {\mathcal B}_{\gamma}(\mathbf{w}_i) = {\mathcal A}(\epsilon)\)
and so from~(\ref{a_eps_est}) and~(\ref{a_set_def}), we have
\begin{equation}\label{q_est}
n^{r}\left(1 - n^{-r\epsilon}\right) \leq \#{\mathcal A}(\epsilon) \leq \sum_{i=1}^{M} \#{\mathcal B}_{\gamma}(\mathbf{w}_i) \leq M \cdot n^{r(\theta+\epsilon)}
\end{equation}
Thus~\(M \geq n^{(1-\theta -\epsilon)r}(1-n^{-r\epsilon})\) and choosing~\(\epsilon = \frac{1}{\sqrt{r}},\) taking~\(r^{th}\) roots and allowing~\(r \rightarrow \infty,\) we get that
\begin{equation}\label{d_av_bd}
\Theta_{\gamma}(G) \geq f(\gamma,d_{av})
\end{equation}
for~\(0 < \gamma < 1-\frac{1}{d_{av}+1}.\)

For~\(1-\frac{1}{d_{av}+1} \leq \gamma \leq 1,\) we use the fact that the expected ball size in~(\ref{eq_one2}) is bounded above by
\begin{equation}\label{eq_one3}
\mathbb{E}\#{\mathcal B}_{\gamma}(\mathbf{v}) \leq \sum_{k=0}^{r} {r \choose k} d^{k}_{av} = (d_{av}+1)^{r}.
\end{equation}
As before we then use the maximality argument to get that~(\ref{d_av_bd}) holds for~\(1-\frac{1}{d_{av}+1} \leq \gamma \leq 1.\)

To see that~\(\Theta_{\gamma}(G) \geq f(\gamma,\Delta),\) we use a similar argument as above along with the estimate that for any deterministic vector~\(\mathbf{u} \in V^{r}\) the ball size
\begin{equation}\label{eq_one4}
\#{\mathcal B}_{\gamma}(\mathbf{u}) \leq \sum_{k=0}^{\gamma r} {r \choose k} \Delta^{r}.
\end{equation}
Finally, a similar argument also shows that~\(\Theta_{\gamma}(G) \geq f(\gamma,L-1)\) and this completes the proof of the lower bound in~(\ref{thet_bds}).\(\qed\)

\emph{Proof of the upper bound in~(\ref{thet_bds})}: Let~\(d = \gamma r\) and let~\({\mathcal C} \subset G(n)\) be a stable set of maximum size in~\(G(r,d)\) and set~\(A(r,d) := \#{\mathcal C}.\) We first derive a recursive estimate involving~\(A(r,d).\) For a vertex~\(v \in G\) let~\({\mathcal C}(v)\) be the set of all vertices in~\({\mathcal C}\) whose last entry is~\(v.\) By the pigeonhole principle, there exists a~\(v_0 \in G\) such that~\(\#{\mathcal C}(v_0) \geq \frac{A(r,d)}{n}.\) We remove the last entry in each vertex of~\({\mathcal C}(v_0)\) and call the resulting set as~\({\mathcal D}(v_0) \subset G(r-1).\) The set~\({\mathcal D}(v_0)\) is also stable in~\(G(r-1,d)\) and this implies that
\[A(r,d) \leq n\cdot A(r-1,d).\] Continuing iteratively we get
\begin{equation}\label{and}
A(r,d) \leq n^{r-d} \cdot A(d,d).
\end{equation}
Taking~\(r^{th}\) roots and using the fact that
\[\sup_{r \geq \frac{1}{\gamma}}(A(d,d))^{\frac{1}{r}} = \sup_{d \geq 1}(A(d,d))^{\gamma d} = (\Theta(G))^{\gamma}, \]
we get the upper bound in~(\ref{thet_bds}).~\(\qed\)

\section*{Acknowledgements}
I thank Professors V. Guruswami, C. R. Subramanian and the referees for crucial comments that led to an improvement of the paper. I also thank IMSc for my fellowships.

\bibliographystyle{plain}

\end{document}